\documentclass[10pt]{article}
\usepackage[below]{placeins}

\usepackage{amsmath}
\usepackage{amsfonts}
\usepackage{amssymb}

\newtheorem{proposition}{Proposition}[section]
\newtheorem{theorem}{Theorem}[section]

\newtheorem{remark}[theorem]{Remark}
\newtheorem{example}[theorem]{Example}

\def\phi{{\varphi}}

\DeclareSymbolFont{AMSb}{U}{msb}{m}{n}
\DeclareMathSymbol{\N}{\mathbin}{AMSb}{"4E}
\DeclareMathSymbol{\Z}{\mathbin}{AMSb}{"5A}
\DeclareMathSymbol{\R}{\mathbin}{AMSb}{"52}
\DeclareMathSymbol{\Q}{\mathbin}{AMSb}{"51}
\DeclareMathSymbol{\I}{\mathbin}{AMSb}{"49}
\DeclareMathSymbol{\C}{\mathbin}{AMSb}{"43}

\begin{document}

\title{Simultaneous preconditioning and symmetrization of non-symmetric linear systems}

\author{ Nassif  Ghoussoub\thanks{Partially supported by a grant
from the Natural Sciences and Engineering Research Council of Canada.  } \quad  and \quad Amir %%@
Moradifam \thanks{Partially supported by a UBC Graduate Fellowship.  }
\\
\small Department of Mathematics,
\small University of British Columbia, \\
\small Vancouver BC Canada V6T 1Z2 \\
\small {\tt nassif@math.ubc.ca} \\
\small {\tt a.moradi@math.ubc.ca}
\\
%\today\\
%\date{January 20, 2005}\\
}
\maketitle

\begin{abstract}
 Motivated by the theory of self-duality  which provides a variational formulation and resolution %%@
for  non self-adjoint partial differential equations \cite{G1, G2}, we propose new templates for  %%@
solving  large non-symmetric linear systems. The method consists of combining a new scheme that  %%@
simultaneously preconditions and symmetrizes the problem,  with various well known iterative methods %%@
for solving linear and symmetric problems. The approach seems to be efficient when dealing with certain  ill-conditioned, and highly non-symmetric systems.  
\end{abstract}

\section{Introduction and main results} 
Many problems in scientific computing lead to systems of linear equations of the form, 
\begin{equation}\label{main}
\hbox{$Ax=b$ where $A\in \R^{n\times n}$ is a nonsingular but sparse matrix, and $b$ is a given vector in %%@
$\R^{n}$,  }
\end{equation} 
%where $A$ is a large sparse matrix,  
and various iterative methods have  been developed for %%@
a fast and efficient resolution of such systems. The Conjugate Gradient Method (CG) which is the oldest and best known of the nonstationary iterative methods, is highly effective in solving symmetric positive definite systems. For indefinite matrices, the minimization feature of CG is no longer an option,  but the Minimum Residual (MINRES) and the Symmetric LQ (SYMMLQ) methods are often computational alternatives for CG, since they are applicable to  systems whose coefficient matrices are symmetric but possibly indefinite.

The case of non-symmetric linear systems is more challenging,  and again methods  such as CGNE, CGNR, GMRES, BiCG, QMR, CGS, and Bi-CGSTAB have been developed %%@
to deal with these situations (see the survey books \cite{Gr} and \cite{S}). 
One approach to deal with the non-symmetric case, consists of reducing the problem to a symmetric one to which one can apply the  above mentioned schemes.  The one that is normally used consists of simply applying CG to the normal 
equations
\begin{equation}\label{normal-equ}
A^{T}Ax=A^{T}b \ \ \hbox{ or} \ \ AA^{T}y=b, \ \ x=A^{T}y.
\end{equation}
It is easy to understand and code this approach, and the CGNE and CGNR methods are based on this  %%@
idea. However, the convergence analysis of these methods depends closely on the {\it condition number}  of the matrix under study. For a general matrix $A$, the condition number  is  
defined as 
\begin{equation}
\kappa (A)=\|A\|\cdot \|A^{-1}\|, 
\end{equation}
and in the case where  $A$ is positive definite and symmetric, the condition number is then equal to  \begin{equation}
\tilde \kappa (A)=\frac{\lambda_{\rm max}(A)}{\lambda_{\rm min}(A)}, 
\end{equation}
 where $\lambda_{\rm min}(A)$ (resp., $\lambda_{\rm max}(A)$) is the smallest %%@
(resp., largest) eigenvalue of $A$). 
% we echo Greenbaum's statement that numerical analysts {\it cringe} at the thought of solving these normal equations because  the {\it condition number} (see below) of the new matrix $A^{T}A$ is the square of the condition number of the original matrix $A$.
%, making  the convergence of CG very slow. 
%This in spite of the futility of using the condition number of $A$ in this comparison since there is no iterative method for non-Hermitian problems whose convergence rate is governed by the condition number of the original  matrix $A$. Note that the {\it condition number} of a general matrix $A$ is  
%defined as 
%\begin{equation}
%\kappa (A)=\|A\|\cdot \|A^{-1}\|, 
%\end{equation}
%and in the case where  $A$ is positive definite and symmetric, the condition number is then equal to  \begin{equation}
%\tilde \kappa (A)=\frac{\lambda_{\rm min}(A)}{\lambda_{\rm max}(A)}, 
%\end{equation}
% where $\lambda_{\rm min}(A)$ (resp., $\lambda_{\rm max}(A)$) is the smallest %%@
%(resp., largest) eigenvalue of $A$). 
The two expressions can be very different for non-symmetric matrices, and these are precisely the systems that seem to be the most pathological from the numerical point of view. Going back to the crudely symmetrized system (\ref{normal-equ}), we echo Greenbaum's statement  \cite{Gr} that numerical analysts {\it cringe} at the thought of solving these normal equations because  the {\it condition number} (see below) of the new matrix $A^{T}A$ is the square of the condition number of the original matrix $A$.
%  and the ones that are least amenable to the known methods. 

In this paper, we shall follow a similar approach that consists of symmetrizing the problem  so as to be able to apply CG, MINRES, or SYMMLQ. However, we argue that for a large class  of non-symmetric, ill-conditionned matrices, it is sometimes beneficial to replace problem (\ref{main}) 
by one of the form  
\begin{equation}\label{gen1-equ-system}
A^{T}MAx=A^{T}Mb,  
\end{equation}
where $M$ is a symmetric and positive definite matrix that can be chosen properly so as to %%@
obtain good convergence behavior for CG when it is applied to the resulting symmetric $A^{T}MA$. This reformulation should not only be seen as a symmetrization, but also as preconditioning procedure.  While it is difficult to %%@
obtain general conditions on $M$ that ensure higher efficiency by minimizing the condition number %%@
$k(A^{T}MA)$, we shall show theoretically and numerically  that by choosing $M$ to be  either the inverse of the symmetric  part of $A$,  or its resolvent, one can get surprisingly good numerical schemes to solve (\ref{main}).

The basis of our approach  originates from the  selfdual variational principle developed in %%@
\cite{G1, G2} to provide a variational formulation and resolution for  non self-adjoint partial %%@
differential equations that do not normally fit in the standard Euler-Lagrangian theory. Applied %%@
to the linear system (\ref{main}), the new principle yields the following procedure. Split %%@
the  matrix $A$ into its  symmetric $A_{a}$ (resp., anti-symmetric part $A_{a}$)
\begin{equation}\label{split1}
A=A_{s}+A_{a},
\end{equation}
where
\begin{equation}\label{split2} 
A_{s}:=\frac{1}{2}(A+A^{T}) \ \ \hbox{and} \ \ A_{a}:=\frac{1}{2}(A-A^{T}).
\end{equation}

\begin{proposition}  {\rm (Selfdual symmetrization)} Assume the matrix $A$ is positive definite, %%@
i.e., for some $\delta>0$, 
\begin{equation} \label{coercive}
\hbox{ $\langle Ax, x\rangle \geq \delta |x|^{2}$ for all $x\in \R^{n}$.}  
\end{equation}
The convex continuous functional 
\begin{equation}
I(x)=\frac{1}{2}\langle Ax,x\rangle +\frac{1}{2}\langle A^{-1}_{s}(b-A_{a}x), b-A_{a}x\rangle %%@
-\langle b, x\rangle
\end{equation}
then attains its minimum at some $\bar x$ in $\R^n$, in such a way that 
\begin{eqnarray}
I(\bar x)&=&\inf\limits_{x\in \R^n}I(x)=0\\
A \bar x&=&b. 
\end{eqnarray}
\end{proposition}
{\bf Symmetrization and preconditioning via selfduality:} Note that the functional $I$ can be written as
\begin{equation}
I(x)=\frac{1}{2}\langle \tilde Ax,x\rangle +\langle A_aA_s^{-1}b-b, x\rangle +\frac{1}{2}\langle A_s^{-1}b, b\rangle,
 \end{equation}
where
\begin{equation}
 \tilde A:=A_{s}-A_{a}A_{s}^{-1}A_{a}=A^{T}A_{s}^{-1}A.
%, \quad \tilde b:=A_aA_s^{-1}b-b$ \, and  $\tilde c:=\frac{1}{2}\langle A_s^{-1}b, b\rangle$.}  
\end{equation}
By writing that $DI(\bar x)=0$, one gets the following equivalent way of solving (\ref{main}).  \\
%\begin{proposition} 

{\em  If both  $A \in %%@
\R^{n\times n}$ and its symmetric part $A_s$ are nonsingular, then $x$ is a solution of the %%@
equation (\ref{main}) if and only if it is a  solution of the linear symmetric equation
\begin{equation}\label{equ-system}
A^{T}A_{s}^{-1}Ax=(A_{s}-A_{a}A_{s}^{-1}A_{a})x=b-A_{a}A_{s}^{-1}b=A^{T}A_{s}^{-1}b.
\end{equation}
}
%\end{proposition}
One can therefore apply to (\ref{equ-system}) all known iterative methods for symmetric systems to %%@
solve the non-symmetric linear system (\ref{main}). 
As mentioned before,  the new equation (\ref{equ-system}) can be seen as  a new symmetrization of problem %%@
(\ref{main}) which also preserves positivity, i.e., $A^{T}A_{s}^{-1}A$ is positive definite if $A$ %%@
is. This will then allow for the use of the Conjugate Gradient Method (CG) for the functional $I$.  %%@
More important and less obvious than the symmetrization effect of $\tilde A$,  is our observation %%@
that for a large class of matrices, the convergence analysis on the system (\ref{equ-system})  is %%@
often more favorable than the original one. The Conjugate Gradient method --which can %%@
now be applied to the symmetrized matrix $\tilde A$--  has the potential of providing an efficient %%@
algorithm for resolving non-symmetric linear systems. We shall call this scheme the {\it Self-Dual %%@
Conjugate Gradient  for Non-symmetric matrices} and we will refer to it as SD-CGN. 

As mentioned above, the convergence analysis of this method depends closely on the condition number $k (\tilde A)$ of $\tilde %%@
A=A^{T}A_{s}^{-1}A$ which in this case is equal to $\tilde k (\tilde A)$. 
%which is defined as 
%\begin{equation}
%k (\tilde A)=\frac{\lambda_{\rm min}(\tilde A)}{\lambda_{\rm max}(\tilde A)}, 
%\end{equation}
% where $\lambda_{\rm min}(\tilde A)$ (resp., $\lambda_{\rm max}(\tilde A)$) is the smallest %%@
%(resp., largest) eigenvalue of $\tilde A$).  
   We observe in section 2.3 
   %on many basic examples below,  
   that even though $ k (\tilde A)$ could be as large as the square of  $k (A_s)$, it is still much smaller that 
    the condition number of the original matrix $\kappa (A)$.
    % and surprisingly we observe that it can improve as the antisymmetric part of the matrix gets larger. 
    In other words, 
   %$\tilde k (\tilde A)$ 
   % of the matrix in %%@ (\ref{equ-system}) 
%   can be much smaller that the condition number of the original matrix $A$ %%@
%(Section 2.3), indicating that 
the inverse $C$ of  $A^{T}A_{s}^{-1}$ can be an  %%@
efficient preconditioning matrix, in spite of the additional cost involved in finding the inverse %%@
of $A_s$. Moreover, %and contrary to other pre-conditioners, 
the efficiency of $C$ seems to surprisingly improve %%@
in many cases as the norm of the anti-symmetric part gets larger (Proposition 2.2). A typical example is when the anti-symmetric matrix $A_a$ is a multiple of the symplectic matrix $J$ (i.e. $JJ^*=-J^2=I$). Consider then a matrix $A_\epsilon= A_s+\frac{1}{\epsilon}J$ which has an arbitrarily large anti-symmetric part. One can show that \begin{equation}
\kappa (\tilde A_\epsilon) \leq \kappa (A_s) +\epsilon^2 \lambda_{\rm max}(A_s)^2, 
\end{equation}
which means that the larger the anti-symmetric part, the more efficient is our proposed selfdual preconditioning. 
Needless to say that this method is of practical interest only when the equation $A_sx=d$ can be solved with less computational effort than the original system, which is not always the case.

Now the relevance of this approach stems from the fact that 
%We cannot formally show that the condition number of system (\ref{main}) is better than %%@
%(\ref{equ-system}), even though as we mention above, we  can observe it on many examples. Still, %%@
conjugate gradient methods for nonsymmetric systems are costly since they require the storage of %%@
previously calculated vectors. It is however worth noting  that  Concus and Golub \cite{CG} and Widlund \cite{W} have also proposed another way to combine CG with a preconditioning using the symmetric  part $A_s$,  %of the coefficient matrix $A$, 
%, as proposed by Concus and Golub \cite{CG} and Widlund %%@\cite{W} 
which does not need this extended storage. Their method has essentially the same %%@
cost per iteration as the preconditioning with the inverse of $A^TA^{-1}_s$ that we propose for SD-CGN and both 
schemes converge to the solution in at most $N$ iterations. \\

{\bf  Iterated preconditioning:} Another way to see the relevance of $A_s$ as a preconditioner, is by noting 
that the convergence of ``simple iteration" 
\begin{equation}\label{dec}
A_{s}x_k=-A_{a}x_{k-1}+b
\end{equation} 
applied to the decomposition of $A$ into its symmetric and anti-symmetric parts,  
requires that the spectral radius $\rho (I-A_s^{-1}A)=\rho (A_s^{-1}A_a) <1$. By multiplying %%@
(\ref{dec}) by $A_s^{-1}$, we see  that this is equivalent to the process of applying simple %%@
iteration  to the original system (\ref{main}) conditioned by $A_s^{-1}$, i.e., to the system %%@
\begin{equation}\label{golub}
A_s^{-1}Ax=A_s^{-1}b.
\end{equation}
On the other hand, ``simple iteration" applied to the decomposition of $\tilde A$ into $A_s$ and  $A_{a}A_{s}^{-1}A_{a}$ is given by%(\ref{equ-system}) 
\begin{equation}
A_{s}x_k=A_{a}A_{s}^{-1}A_{a}x_{k-1}+b-A_{a}A_{s}^{-1}b.
\end{equation}
Its convergence   is controlled by $\rho (I-A_s^{-1}\tilde  A)=\rho ((A_s^{-1}A_a)^2)=\rho(A_s^{-1}A_a)^2$ which is strictly less than $\rho(A_s^{-1}A_a)$, i.e., an improvement when  the latter is strictly less than one,  which  the mode in which we have convergence. In other words, the linear system (\ref{equ-system}) can still be preconditioned one more time as follows:\\

 {\em If both  $A %%@
\in \R^{n\times n}$ and its symmetric part $A_s$ are nonsingular, then $x$ is a solution of the %%@
equation (\ref{main}) if and only if it is a  
solution of the linear symmetric equation
\begin{equation}\label{equ-system.2}
\bar %%@
Ax:=A_s^{-1}A^{T}A_{s}^{-1}Ax=[I-(A_{s}^{-1}A_{a})^2]x=(I-A_s^{-1}A_{a})A_{s}^{-1}b=A_s^{-1}A^{T}A
_{s}^{-1}b.
\end{equation}
} 
Note however that with this last formulation, one has to deal with the potential loss of  positivity for the matrix $\tilde A$.\\ 
%It is worth noting here that the benefits of this iterated pre-conditioning, may unfortunately be offset by the potential loss of  positivity for the matrix $\tilde A$.\\ 
 %\end{proposition}
%This scheme will be denoted PSD-CGN, and we shall see in section 2 that it is more efficient than %%@
%the preconditioner  of Concus-Golub, albeit with simple iteration or with CG.  

{\bf Anti-symmetry in transport problems:} 
Numerical experiments on standard linear ODEs (Example 3.1) and  PDEs (Example 3.2), show the efficiency of
SD-CGN  for non-selfadjoint equations. %%@
Roughly speeking, discretization of differential equations normally leads to a symmetric component %%@
coming from the Laplace operator, while the discretization of the  non-self-adjoint part leads to %%@
the anti-symmetric part of the coefficient matrix. As such, the symmetric part of the matrix is of %%@
order $O(\frac{1}{h^2})$, while the anti-symmetric part is of order $O(\frac{1}{h})$, where $h$ is %%@
the step size. The coefficient matrix $A$ in the original system (\ref{main}) is therefore an %%@
$O(h)$ perturbation of its symmetric part. However, for the new system (\ref{equ-system}) we have %%@
roughly 
\begin{equation}
\tilde %%@
A=A_{s}-A_{a}A_{s}^{-1}A_{a}=O(\frac{1}{h^2})-O(\frac{1}{h})O(h^2)O(\frac{1}{h})=O(\frac{1}{h^2})-%%@
O(1),
\end{equation}
making the matrix $\tilde A$ an $O(1)$ perturbation of $A_{s}$, and therefore a matrix of the form  $A_{s}+\alpha I$  becomes a natural candidate to precondition the new system (\ref{equ-system}).   \\

{\bf Resolvents of $A_s$ as preconditioners:} One may therefore consider preconditioned equations of the form 
%in full generality 
%We shall try in the sequel to compare the effect of the selfdual symmetrization and %%@
%preconditioning procedure (\ref{equ-system}) to others in the literature including the standard %%@
%one 
%in (\ref{normal-equ}). We shall see that it often leads to more efficient algorithms, but we are %%@
%however unable to classify the class of matrices $A$ for which it is optimal.  More generally, one %%@
%can always embed both procedures in one of the form  
%\begin{equation}\label{gen1-equ-system}
$A^{T}MAx=A^{T}Mb$,   
%\end{equation}
where $M$ is of the form
%can 
% is a symmetric and positive definite matrix that needs to be chosen properly so as to %%@
%obtain good convergence behavior for CG when it is applied to $A^{T}MA$. While it is difficult to %%@
%obtain general conditions on $M$ that ensure higher efficiency by minimizing the condition number %%@
%$k(A^{T}MA)$, we shall consider numerical experiments involving matrices of the form   
\begin{equation}
\hbox{$M_{\alpha}=\big(\alpha A_s+(1-\alpha) I\big)^{-1}$\quad  or \quad $N_{\beta}=\beta %%@
A_s^{-1}+(1-\beta)I$,}
\end{equation}
for some $0\leq \alpha, \beta \in \R$, and where $I$ is the unit matrix.  

Note that we obviously recover (\ref{normal-equ}) when $\alpha =0$, and  (\ref{equ-system}) when %%@
$\alpha =1$. 
As $\alpha \rightarrow 0$ the matrix $\alpha A_s+(1-\alpha) I$ becomes easier to invert, but the %%@
matrix 
\begin{equation}
A_{1, \alpha}=A^{T}(\alpha A_s+(1-\alpha) I)^{-1}A
\end{equation}
may become more ill conditioned, eventually leading (for $\alpha=0$) to 
$A^{T}Ax=A^{T}b$. There is therefore a trade-off between the efficiency of CG for the system %%@
(\ref{gen1-equ-system}) and the condition number of the inner matrix $\alpha A_s+(1-\alpha) I$, %%@
and so by an appropriate choice of the parameter $\alpha$ we may minimize the cost of finding a %%@
solution for the system (\ref{main}). In the case where $A_s$ is positive definite, one can choose %%@
--and it is sometimes preferable as shown in example (3.4)-- $\alpha>1$, as long as  $\alpha %%@
<\frac{1}{1-\lambda^s_{\rm min}}$,  where $\lambda^{s}_{\rm min}$ is the smallest eigenvalue of %%@
$A_s$. Moreover, in the case where the matrix $A$ is not positive definite or if its symmetric %%@
part is not invertible,  one may take $\alpha$ small enough, so that the matrix $M_{\alpha}$ (and %%@
hence $A_{1, \alpha}$) becomes positive definite, and therefore making CG applicable (See example 3.4). 
Similarly,  the matrix $N_{\beta}=\beta A_s^{-1}+(1-\beta)I$ provides another choice for the %%@
matrix $M$ in (\ref{gen1-equ-system}), for $\beta <\frac{\lambda^{s}_{\rm max}}{\lambda^{s}_{\rm %%@
max}-1}$ where $\lambda^{s}_{\rm max}$ is the largest eigenvalue of $A_s$. Again we may choose %%@
$\alpha$ close to zero to make the matrix $N_{\beta}$  
positive definite. As we will see in the last section, appropriate choices of $\beta$, can lead to %%@
better convergence of CG for equation (\ref{gen1-equ-system}).

One can also combine both effects by considering matrices of the form
\begin{equation}
\hbox{$L_{\alpha, \beta}=\big(\alpha A_s+(1-\alpha) I\big)^{-1} + \beta I$,}
\end{equation}
as is done in example (3.4). 

We also note that the matrices $M'_{\alpha}:=(\alpha A'_s+(1-\alpha)I)^{-1}$ and %%@
$N'_{\beta}:=\beta(A'_s)^{-1}+(1-\beta)I$ can be  other options for the matrix $M$, where $A'_s$ %%@
is a suitable approximation of $A_s$, chosen is such a way that $M'_{\alpha}q$ and $N'_{\beta}q$ %%@
can be relatively easier to compute for any given vector $q$. 

Finally, we observe that the above reasoning applies to any decomposition $A=B+C$ of the %%@
non-singular matrix $A \in \R^{n\times n}$, where  $B$ and $(B-C)$ are both invertible. In this %%@
case,  $B(B-C)^{-1}$ can be a preconditioner for the equation (\ref{main}). Indeed, since 
$B-CB^{-1}C=(B-C)B^{-1}A$,
$x$ is a solution of   (\ref{main}) if and only of it is a solution of the system 
\begin{equation}\label{gen-equ-system}
(B-C)B^{-1}Ax=(B-CB^{-1}C)x=b-CB^{-1}b.
\end{equation}
In the next section, we shall describe a general framework based on the ideas explained above for the %%@
use of  iterative  methods for solving non-symmetric linear systems. In section 3 we present %%@
various numerical experiments to test the  effectiveness of the proposed methods.

\section{Selfdual methods for non-symmetric systems}

By {\it selfdual methods} we mean the ones that  consist of first associating to problem %%@
(\ref{main}) 
the equivalent system (\ref{gen1-equ-system}) with appropriate choices of $M$, then exploiting the %%@
symmetry of the  new system by using the various existing  iterative methods for symmetric systems %%@
such as CG, MINRES, and SYMMLQ, leading eventually to the solution of the original problem %%@
(\ref{main}).
 In the case where the matrix $M$ is positive definite, one can then use CG on the equivalent %%@
system (\ref{gen1-equ-system}). This scheme  (SD-CGN) 
%(Self-Dual Conjugate Gradient %%@ for  Nonsymmetric linear equations) and
 is illustrated in Table (1) below, in the case where the %%@
matrix $M$ is chosen to be the inverse of the symmetric part of $A$. If $M$ is not positive %%@
definite, then one can  use MINRES (or SYMMLQ) to solve the system (\ref{equ-system}). We will %%@
then refer to them as SD-MINRESN (i.e., Self-Dual  MINRES for Nonsymmetric linear equations). 

\subsection{Exact methods}  In each iteration of CG, MINRES, or SYMMLQ, one needs to compute $Mq$ %%@
for certain vectors $q$.
Since selfdual methods call for a conditioning matrix $M$ that involves inverting another one, %%@
the computation of $Mq$ can therefore be costly, and therefore not necessarily efficient for all %%@
linear equations.  But as we will see in section 3, $M$ can sometimes be chosen  so that computing %%@
$Mq$ is much easier than solving the original 
equation itself. This is the case for example when the symmetric part is either diagonal or %%@
tri-diagonal, or when we are dealing with several linear systems all having the same symmetric %%@
part, but with different anti-symmetric components. Moreover, one need not find the whole matrix %%@
$M$, in order to compute $Mq$. The following scheme illustrates the exact SD-CGN method applied in %%@
the case where the coefficient matrix $A$ in (\ref{main}) is positive definite, and when %%@
$A^{T}(A_{s})^{-1}Aq$ can be computed exactly for any given vector $q$.

\begin{table}[ht]
\begin{center}
\begin{tabular}{|l|}\hline \\ 
Given an initial guess $x_{0}$,\\ 
Solve $A_{s}y=b$\\ 
Compute $\overline{b}=b-A_{a}y$.\\ 
Solve $A_{s}y_{0}=A_{a}x_{0}$\\ 
Compute  $r_{0}=\overline{b}-A_{s}x_{0}+A_{a}y_{0}$  and set $p_{0}=r_{0}$.\\ 
For k=1,2, . . . ,\\ 
Solve $A_{s}z=A_{a}p_{k-1}$\\ 
Compute $w=A_{s}p_{k-1}-A_{a}z$ .\\  
Set $x_{k}=x_{k-1}+\alpha_{k-1}p_{k-1}$, where %%@
$\alpha_{k-1}=\frac{<r_{k-1},r_{k-1}>}{<p_{k-1},w>}$ .\\ 
Cpmpute $r_{k}=r_{k-1}-\alpha_{k-1}w$.\\ 
Set $p_{k}=r_{k}+b_{k-1}p_{k-1}$,  where $b_{k-1}=\frac{<r_{k},r_{k}>}{<r_{k-1},r_{k-1}>}$ .\\ 
Check convergence; continue if necessary. \\ \\
\hline
\end{tabular}
\end{center}
\caption{GCGN}
\end{table}
\FloatBarrier

In the case where $A$ is not positive definite, or when it is preferable to choose a non-positive %%@
definite conditioning matrix $M$, then one can apply MINRES or SYMMLQ   to the equivalent system %%@
(\ref{gen1-equ-system}).  These schemes will be then called  SD-MINRESN and SD-SYMMLQN %%@
respectively.

\subsection{Inexact Methods}
The SD-CGN, SD-MINRESN and SD-SYMMLQN are of practical interest when for example, the equation 
\begin{equation}\label{sym}
A_{s}x=q
\end{equation}
can be solved with less computational effort than the original equation (\ref{main}). Actually, %%@
one can use CG, MINRES, or SYMMLQ to solve (\ref{sym}) in every iteration of SD-CGN, SD-MINRESN, %%@
or SD-SYMMLQN. But since each sub-iteration may lead to an error in the computation of %%@
(\ref{sym}), one needs to control such errors, in order for the method to lead to a solution of %%@
the system (\ref{main}) with the desired tolerance. This leads to the Inexact SD-CGN, SD-MINRESN %%@
and SD-SYMMLQN methods (denoted below by ISD-CGN, ISD-MINRESN and ISD-SYMMLQN respectively). 

The following proposition --which is a direct consequence of Theorem 4.4.3 in \cite{Gr}-- shows %%@
that if we solve the inner equations (\ref{sym}) ``accurately enough" then ISD-CGN and ISD-MINRESN %%@
can be used to solve  (\ref{main}) with a pre-determined accuracy. Indeed, given 
$\epsilon>0$, we assume that in each iteration of ISD-CGN or  ISD-MINRESN,  we can solve the inner %%@
equation --corresponding to $A_{s}$--  accurately enough in such a way  that 
\begin{equation}
\|(A_{s}-A_{a}A_{s}^{-1}A_{a})p-(A_{s}p-A_{a}y)\|=\|A_{a}A_{s}^{-1}A_{a}p-A_{a}y\|<\epsilon,
\end{equation}
where $y$ is the (inexact) solution of the equation
\begin{equation}\label{sym.2} 
A_{s}y=A_{a}p.
\end{equation}
In other words, we assume CG and MINRES are implemented on (\ref{sym.2}) in a finite precision %%@
arithmetic with machine precision $\epsilon$. Set 
\begin{equation}
\epsilon_{0}:=2(n+4)\epsilon, \ \ \ \ \epsilon_{1}:=2(7+n\frac{\|\, |A_{s}-A_{a}A_{s}^{-1}A_{a}|\, \| %%@
|}{\|A_{s}-A_{a}A_{s}^{-1}A_{a}\|})
\epsilon, 
\end{equation}
where  $|D|$  denotes the matrix  whose terms are the absolute values of the corresponding terms in the matrix $D$. 
 Let $\lambda_{1}\leq...\leq \lambda_{n}$ be the eigenvalues of $(A_{s}-A_{a}A_{s}^{-1}A_{a})$ and %%@
let $T_{k+1,k}$ be the $(k+1)\times k$ tridiagonal matrix generated by a finite precision Lanczos %%@
computation. Suppose that there exists a symmetric tridiagonal matrix  $T$, with $T_{k+1,k}$ as %%@
its upper left $(k+1)\times k$ block, whose eigenvalues all lie in the intervals 
\begin{equation}
S=\cup^{n}_{i=1}[\lambda_{i}-\delta,\lambda_{i}+\delta],
\end{equation} 
where none of the intervals contain the origin. let $d$ denote the distance from the origin to the %%@
set $S$, and let $p_k$ denote a polynomial of degree $k$. 

\begin{proposition}
The ISD-MINRESN residual $r^{IM}_{k}$ then satisfies 
\begin{equation}
 \frac{||r^{IM}_{k}||}{||r_{0}||}\leq \sqrt{(1+2\epsilon_{0})(k+1)}\ \ \min_{p_{k}}\ \ %%@
\max_{z=S}|p_{k}(z)|+2\sqrt{k}(\frac{\lambda_{n}}{d})
\epsilon_{1}.
\end{equation}
If $A$ is positive definite, then the ISD-CGN residual $r^{IC}$ satisfies 
\begin{equation}
 \frac{||r^{IC}_{k}||}{||r_{0}||}\leq \sqrt{(1+2\epsilon_{0})(\lambda_{n}+\delta)/d}\ \ %%@
\min_{p_{k}}\ \ \max_{z=S}|p_{k}(z)|+\sqrt{k}(\frac{\lambda_{n}}{d})
\epsilon_{1}.
\end{equation}
\end{proposition}
It is shown by Greenbaum \cite{G1} that $T_{k+1,k}$ can be extended to a larger symmetric %%@
tridiagonal matrix $T$ whose eigenvalues all lie in tiny intervals about the eigenvalues of 
$(A_{s}-A_{a}A_{s}^{-1}A_{a})$. Hence the above proposition guarantees that if we solve the inner %%@
equations accurate enough, then ISD-CGN and ISD-MINRESN converges to the solution of the system %%@
\ref{main} with the desired relative residual (see the last section for numerical experiments).

\subsection{Preconditioning } 

As mentioned in the introduction, the convergence of iterative methods depends heavily on the %%@
spectral properties of the coefficient matrix. Preconditioning techniques attempt to transform the %%@
linear system (\ref{main}) into an equivalent one of the form $C^{-1}Ax=C^{-1}b$,  in such a way  %%@
that it has the same solution, but hopefully with more favorable spectral properties. As such the %%@
reformulation of (1)  as 
\begin{equation}
A^{T}A_{s}^{-1}Ax=A^{T}A_{s}^{-1}b, 
\end{equation}
can be seen as a preconditioning procedure with $C$ being the inverse of $A^TA_s^{-1}$. The %%@
spectral radius, and more importantly the condition number of the coefficient matrix in linear %%@
systems,  are crucial  parameters for the convergence of iterative methods. The following %%@
simple proposition gives upper bounds on the condition number of  $\tilde A=A^{T}A_{s}^{-1}A$.
 
\begin{proposition}Assume $A$ is an invertible positive definite matrix, then  
\begin{equation}\label{cond-est}
\kappa (\tilde A) \leq \min \{\kappa_1, \kappa_2\}, 
\end{equation}
where
\begin{equation}
\hbox{$\kappa_1:=\kappa (A_{s})+\frac{\|A_a\|^2}{\lambda_{\rm min}(A_{s})^2}$ \quad and \quad $  \kappa_2:=\kappa(A_{s})\kappa(-A_a^2)+\frac{\lambda_{\rm max}(A_{s})^2}{\lambda_{\rm min} (-A_a^2)}$.}
\end{equation}
\end{proposition}
{\bf Proof:} We have
\begin{eqnarray*}
\lambda_{min}(\tilde A)=\lambda_{min}(A_{s}-A_{a}A_s^{-1}A_{a})\geq \lambda_{min}(A_{s}).
\end{eqnarray*}
We also have 
\begin{eqnarray*}
\lambda_{max}(\tilde A)&=&\sup_{x\neq 0} \frac{x^{t}\tilde A x}{|x|^2}=\sup_{x\neq 0} %%@
\frac{x^{t}(A_{s}-A_{a}A_s^{-1}A_{a}) x}{|x|^2}\\
&\leq& \lambda_{max}(A_{s})+\frac{||A_{a}||^{2}}{\lambda_{min}(A_{s})}.
\end{eqnarray*}
Since $\kappa (\tilde A)=\frac{\lambda_{max}(\tilde A)}{\lambda_{min}(\tilde A)}$, it follows that $\kappa (\tilde A) \leq \kappa_1$.

To obtain the second estimate, observe that
\begin{eqnarray*}
\lambda_{min}(\tilde A)&=&\lambda_{min}(A_{s}-A_{a}A_s^{-1}A_{a})>\lambda_{min}(-A_{a}A_s^{-1}A_{a})\\
&=&\inf_{x\neq 0}\frac{-x^{T}A_{a}A_s^{-1}A_{a}x}{x^{T}x}\\
&=&\inf_{x\neq 0}\{\frac{(A_{a}x)^{T}A_s^{-1}(A_{a}x)}{(A_{a}x)^{T}(A_{a}x)}\times \frac{(A_{a}x)^{T}(A_{a}x)}    %%@
{x^{T}x}\}\\
&\geq& \inf_{x\neq 0}\frac{(A_{a}x)^{T}A_s^{-1}(A_{a}x)}{(A_{a}x)^{T}(A_{a}x)}\times \inf_{x\neq 0} %%@
\frac{x^{T}(A_{a})^{T}(A_{a})x} {x^{T}x}\\
&=&\frac{1}{\lambda_{max}(A_{s})}\times %%@
\lambda_{min}((A_{a})^{T}A_{a})\\
&=&\frac{1}{\lambda_{max}(A_{s})}\times %%@
\lambda_{min}(-A_{a}^2)
%=\frac{||(A_{a})^{T}A_{a}||}{\lambda_{max}(A_{s})}=\frac{||A_{a}^2||}{\lambda_{max}(A_{%%@
%s})}=\frac{||A_{a}||^2}{\lambda_{max}(A_{s})}.
\end{eqnarray*}
With the same estimate for $\lambda_{max}(\tilde A)$ we get $\kappa (\tilde A) \leq \kappa_2$.

%Why $||A_{a}^2||=||A_{a}||^2$?\\
%Proof.
%\begin{eqnarray*}
%||A_{a}^2||=\inf_{x\neq 0}\frac{|A_{a}(A_{a}x)|}{|x|}\geq\inf_{x\neq 0}\frac{|A_{a}(A_{a}x)|}{|A_{a}x|}\inf_{x\neq %%@
%0}\frac{|A_{a}x|}{|x|}=||A_{a}||^2,
%\end{eqnarray*}
%and $||A_{a}^2||\leq||A_{a}||^2$ is trivial.

\begin{remark} \rm Inequality (\ref{cond-est}) shows that SD-CGN and SD-MINRES can be very efficient schemes for a large class of ill conditioned non-symmetric matrices, even those that are almost singular and with arbitrary large condition numbers.  It suffices that either $\kappa_1$ or $\kappa_2$ be small. Indeed, 
 \begin{itemize}
\item The inequality $\kappa (\tilde A) \leq \kappa_1$ shows that the condition number $\kappa (\tilde A)$ is reasonable as long as the anti-symmetric part $A_a$ is not too large. On the other hand, even if $\|A_a\|$ is of the order of $\lambda_{\rm max}(A_s)$, and $\kappa (\tilde A)$ is then as large as $\kappa (A_s)^2$,  it may still be an improved situation, since this can happen for cases when $\kappa (A)$ is exceedingly large. This can be seen in example 2.2 below.

\item  The inequality $\kappa (\tilde A) \leq \kappa_2$ is even more interesting especially  in situations when $\lambda_{\rm min}(-A_a^2)$ is arbitrarily large while remaining of the same order as  $||A_{a}||^2$. This means that  $\kappa (\tilde A)$ can remain of the same order as  $\kappa (A_s)$ regardless how large is $A_a$.   

A typical example is when the anti-symmetric matrix $A_a$ is a multiple of the symplectic matrix $J$ (i.e. $JJ^*=-J^2=I$). Consider then a matrix $A_\epsilon= A_s+\frac{1}{\epsilon}J$ which has an arbitrarily large anti-symmetric part. By using that $\kappa (\tilde A) \leq \kappa_2$, one gets
\begin{equation}
\kappa (\tilde A_\epsilon) \leq \kappa (A_s) +\epsilon^2 \lambda_{\rm max}(A_s)^2.
\end{equation}
  
\end{itemize}   
\end{remark}

Here are other examples where the larger  the condition number of $A$ is, the more efficient is the proposed selfdual preconditioning. 
\begin{example} \label{2-2}\rm
Consider the matrix  
\begin{equation} 
A_\epsilon=\left[
\begin{tabular}{cl}
$1$ & $ -1$ \\
$1$ & $ -1+\epsilon$  
\end{tabular}
 \right] 
\end{equation}
which is a typical example of an ill-conditioned non-symmetric matrix. One can actually show that %%@
$\kappa(A_\epsilon)=O(\frac{1}{\epsilon}) \rightarrow \infty$ as $\epsilon \rightarrow 0$ with %%@
respect to any norm. However, the condition number of the associated selfdual coefficient matrix
\begin{equation*}
\tilde A_\epsilon=A_{s}-A_{a}(A_{s})^{-1}A_{a}=\left[
\begin{tabular}{cc}
$\frac{\epsilon}{\epsilon-1}$ & $0$ \\
$0$ & $\epsilon$ \\
\end{tabular}
 \right]
\end{equation*}
is $\kappa(\tilde A_\epsilon)=\frac{1}{1-\varepsilon}$, and therefore goes to $1$ as $\varepsilon %%@
\rightarrow 0$. Note also that the condition number of the symmetric part of $A_\epsilon$ 
goes to one as $\epsilon \rightarrow 0$. In other words,  the more ill-conditioned problem %%@
$(\ref{main})$ is,  the more efficient the selfdual conditioned system (\ref{equ-system}) is. 

We also observe that $\kappa(A_{s}^{-1}A)$ goes to $\infty$ as $\epsilon$ goes to zero, which %%@
means that besides making the problem symmetric, our proposed conditioned matrix %%@
$A^{T}A_{s}^{-1}A$ has a much smaller  condition number than the matrix $A_{s}^{-1}A$, which %%@
uses $A_s$ as a preconditioner.

 Similarly,  consider the non-symmetric linear system with coefficient matrix 
\begin{equation} 
A_\epsilon=\left[
\begin{tabular}{cc}
$1$ & $-1+\epsilon$\\
$1$ & $-1$\\
\end{tabular}
 \right].
\end{equation}
\end{example}
As $\epsilon \rightarrow 0$, the matrix becomes again more and more ill-conditioned, while the %%@
condition number of its symmetric part converges to one. Observe now that the condition number of %%@
$\tilde A_\epsilon$ also converges to $1$ as $\epsilon$ goes to zero.  This example shows that self-doual preconditioning can also be very efficient for non-positive definite problems. 

\section{Numerical Experiments}
In this section we present some  numerical examples to illustrate the proposed schemes and to %%@
compare them to other known iterative methods for non-symmetric linear systems. 
Our  experiments have been carried out on Matlab (7.0.1.24704 (R14) Service Pack 1).  In all cases %%@
the iteration was started with $x_{0}=0$.

\begin{example} Consider the ordinary differential equation
\begin{equation}\label{ode2}
-\epsilon y''+y'=f(x), \ \ \hbox{on} \ \ [0,1], \ \ y(0)=y(1)=0. 
\end{equation} 
By discretizing this equation with stepsize  $1/65$ and by using backward difference for the 
first order term,  one obtains a nonsymmetric system of linear equations with 64 unknowns. We %%@
present in Table 2 below, the number of iterations needed for various decreasing values of the %%@
residual $\epsilon$. We use ESD-CGN and ISD-CGN (with relative residual $ 10^{-7}$ for the %%@
solutions of the inner equations). We then compare them to the known  methods CGNE, BiCG, QMR, %%@
CGS, and BiCGSTAB for solving non-symmetric linear systems. We also test preconditioned version %%@
of these methods by using the symmetric part of the corresponding matrix as a preconditioner. 

\begin{table}[ht] \label{sol1-ode2} \caption{Number of iterations to find a solution with relative %%@
residual $10^{-6}$ for equation (\ref{ode2}). $f(x)$ is chosen so that $y=x\sin(\pi x)$ is a %%@
solution.}
\begin{center}
\begin{tabular}{|c|c|c|c|c|c|c|c|}\hline
N=64& $\epsilon=10^{-2}$ %%@
&$\epsilon=10^{-3}$&$\epsilon=10^{-4}$&$\epsilon=10^{-6}$&$\epsilon=10^{-10}$&$\epsilon=10^{-16}$ %%@
\\ 
\hline
ESD-CGN& 22& 8 &5 & 4 &3&2 \\ \hline
ISD-CGN($10^{-7}$)& 24 &  9& 6& 4 &3&2 \\ \hline
GCNE& 88& 64 &64 & 64 &64&64 \\ \hline
QMR& 114& $>1000$ &$>1000$ & $>1000$&$>1000$&$>1000$ \\ \hline
PQMR& 34& 51 &50 & 52 &52&52 \\ \hline
BiCGSTAB& 63.5& 78.5 &92.5 & 98.5 &100.5&103.5 \\ \hline
PBiCGSTAB& 26.5& 46.5 &50.5 & 50 &51.5&51.5 \\ \hline
BiCG& 125& $>1000$ &$>1000$ & $>1000$&$>1000$&$>1000$ \\ \hline
PBiCG& 31& 44 &50& 50 &52&52 \\ \hline
CGS& $>1000$& $>1000$ &$>1000$ & $>1000$&$>1000$&$>1000$ \\ \hline
PCGS& 27& 51 &46& 46 &46&48 \\ \hline
\end{tabular}
\end{center}
\end{table}
\FloatBarrier

\begin{table}[ht] \label{sol2-ode2} \caption{Number of iterations to find a solution with relative %%@
residual $10^{-6}$ for equation (\ref{ode2}). $f(x)$ is chosen so that $y=\frac{x(1-x)}{\cos %%@
(x)}$ is a solution, while the  stepsize used is $1/129$.}
\begin{center}
\begin{tabular}{|c|c|c|c|c|c|c|c|}\hline
N=128& $\epsilon=10^{-2}$ %%@
&$\epsilon=10^{-3}$&$\epsilon=10^{-4}$&$\epsilon=10^{-6}$&$\epsilon=10^{-10}$&$\epsilon=10^{-16}$ %%@
\\ 
\hline
ESD-CGN& 37& 11 &6 & 4 &3&2 \\ \hline
ISD-CGN($10^{-7}$)& 38 &  12& 7& 4 &3&2 \\ \hline
GCNE& 266& 140 &128 & 128 &128&128 \\ \hline
QMR& $>1000$& $>1000$ &$>1000$ & $>1000$&$>1000$&$>1000$ \\ \hline
PQMR& 40& 77 &87 & 92 &90&85 \\ \hline
BiCGSTAB& 136.5& 167.5 &241 & 226.5 &233.5&237.5 \\ \hline
PBiCGSTAB& 35.5& 87.5 &106.5 & 109 &110.5&110.5 \\ \hline
BiCG& $>1000$& $>1000$ &$>1000$ & $>1000$&$>1000$&$>1000$ \\ \hline
PBiCG& 37& 76 &84& 89 &85&91 \\ \hline
CGS& $>1000$& $>1000$ &$>1000$ & $>1000$&$>1000$&$>1000$ \\ \hline
PCGS& 34& 80 &96& 91 &94&90 \\ \hline
\end{tabular}
\end{center}
\end{table}
\FloatBarrier
As we see in Tables 2 and  and 3,  a phenomenon similar to Example \ref{2-2} is occuring. As the %%@
problem gets harder ($\epsilon$ smaller), SD-CGN becomes more 
efficient. These results can be compared with the number of iterations that the HSS iteration %%@
method needs to solve equation (\ref{ode2}) (Tables 3,4, and 5 in \cite{BGN}).
\end{example}

\begin{example} Consider the partial differential equation
\begin{equation}\label{pde1}
 -\Delta u+a(x,y) \frac{\partial u}{\partial x}=f(x,y), \ \ 0\leq x\leq 1, \ \ 0\leq y \leq 1,
\end{equation}
with Dirichlet boundary condition. 
\end{example}
The number of iterations that ESD-CGN and ISD-CGN needed to find a solution with relative residual %%@
$10^{-6}$,  are presented in Table 4 below for different coefficients $a(x,y)$.

\begin{table}[ht] \label{sol-pde1} \caption{Number of iterations (I) for the backward scheme method to %%@
find a solution with relative residual $10^{-6}$ for equation (\ref{pde1}) (Example 3.2)}
\begin{center}
\begin{tabular}{|c|c|c|c|c|c|c|c|}\hline
a(x,y)& N &I (ESD-CGN)&I (ISD-CGN)&Solution \\ 
\hline
100& 49& 18 &18 & random \\ \hline
100& 225 &  40& 37& random \\ \hline
100& 961& 44 &46 & random \\ \hline
100& 961&52& 51 &$\sin \pi x \sin \pi y .\exp((x/2+y)^3)$ \\ \hline
1000& 49& 10 &10 &random \\ \hline
1000& 225& 31 &31 & random \\ \hline
1000& 961& 36 &37 & random \\ \hline
1000& 961& 31 &39 & $\sin \pi x \sin \pi y .\exp((x/2+y)^3)$ \\ \hline
$10^6$& 49& 4 &4& random \\ \hline
$10^6$& 225& 6 &6& random \\ \hline
$10^6$& 961& 6 &6& random \\ \hline
$10^6$& 961& 6 &6& $\sin \pi x \sin \pi y .\exp((x/2+y)^3)$ \\ \hline
$10^{16}$& 961& 2 &2& $\sin \pi x \sin \pi y .\exp((x/2+y)^3)$ \\ \hline
\end{tabular}
\end{center}
\end{table}
\FloatBarrier
\begin{table}[ht] \label{sol-pde1} \caption{Number of iterations (I) for the centered difference %%@
scheme method for equation (\ref{pde1}) (Example 3.2)}
\begin{center}
\begin{tabular}{|c|c|c|c|c|c|c|c|}\hline
a(x,y)& N &I (ESD-CGN)&Solution& Relative Residoual \\ 
\hline
1& 49& 21 & random & $6.71\times 10^{-6}$ \\ \hline
1& 225& 73 & random & $9.95 \times 10^{-6}$ \\ \hline
1& 961& 91 & random & $8.09\times 10^{-6}$ \\ \hline
1& 961&72 &$\sin \pi x \sin \pi y .\exp((x/2+y)^3)$& $9.70\times 10^{-6}$\\ \hline
10& 49& 18 & random & $9.97\times 10^{-6}$ \\ \hline
10& 225& 65 & random & $5.90\times 10^{-6}$ \\ \hline
10& 961& 78 & random & $8.95\times 10^{-6}$ \\ \hline
10& 961&65 &$\sin \pi x \sin \pi y .\exp((x/2+y)^3)$& $7.78\times 10^{-6}$ \\ \hline
100& 49& 31 & random & $6.07\times 10^{-6}$ \\ \hline
100& 225& 42 & random & $5.20\times 10^{-6}$\\ \hline
100& 961& 43 & random & $5.03\times 10^{-6}$\\ \hline
100& 961&38 &$\sin \pi x \sin \pi y .\exp((x/2+y)^3)$& $4.69\times 10^{-6}$ \\ \hline
1000& 49& 65 & random & $4.54\times 10^{-6}$ \\ \hline
1000& 225& 130 & random & $8.66\times 10^{-6}$\\ \hline
1000& 961& 140 & random & $2.12\times 10^{-6}$\\ \hline
100& 961&150&$\sin \pi x \sin \pi y .\exp((x/2+y)^3)$& $5.98\times 10^{-6}$ \\ \hline
\end{tabular}
\end{center}
\end{table}
\FloatBarrier
Table 4 and 5 can be compared with Table 1 in \cite{W}, where Widlund had tested his Lanczos %%@
method for non-symmetric linear systems. Comparing Table 5 with Table 1 in \cite{W} we see that %%@
for small $a(x,y)$ (1 and 10) Widlund's method is more efficient than SD-CGN, but for large values %%@
of $a$, SD-CGN turns out to be more efficient than Widlund's Lanczos method.  
\begin{remark} As we see in Tables 2,3, and 4, the number of iterations for ESD-CGN and ISD-CGN %%@
(with relative residual $10^{-7}$ for the solutions of the inner equations) are almost the same %%@
One might choose dynamic  relative residuals for the solutions of inner equations to decrease the %%@
average cost per iterations of ISD-CGN. It is interesting to figure out whether there is a %%@
procedure to determine the accuracy of solutions for the inner equations to minimize the total %%@
cost of finding a solution.   

\end{remark}

\begin{example} Consider the partial differential equation
%\begin{equation}\label{pde2}
%-\Delta u+[\partial_{x}(20\exp(3.5(x^2+y^2))u)+20\exp(3.5(x^2+y^2))\partial_{x} u]/2=f(x), \ \ \ \ %%@
%\hbox{on} \ \ [0,1]\times [0,1],
%\end{equation}
\begin{equation}\label{pde2}
-\Delta u+ 10\frac{\partial (\exp(3.5(x^2+y^2)u)}{\partial x}+10\exp(3.5(x^2+y^2))\frac{\partial u}{\partial x}=f(x), \ \ \ \ %%@
\hbox{on} \ \ [0,1]\times [0,1],
\end{equation}
with Dirichlet boundary condition,  and choose $f$ so that $\sin(\pi x)\sin(\pi y)\exp((x/2+y)^3)$ %%@
is  the solution of the equation. We take the stepsize $h=1/31$ which leads to a linear system %%@
$Ax=b$ with 900 unknowns. Table 5 includes the number of iterations which CG needs to converge to %%@
a solution with relative residual $10^{-6}$ when applied to the preconditioned matrix 
\begin{equation}\label{pde2-system}
A^{T}(\alpha A_s^{-1}+(1-\alpha) I)A.
\end{equation}
Table 5 can be compared with Table 1 in \cite{W}, where Widlund has presented the number of %%@
iterations needed to solve equation (\ref{pde2}).
\end{example}

\begin{table}[ht] \label{sol1-pde2} \caption{Number of iterations for a solution with relative %%@
residual $10^{-6}$ for  example 3.3 when SD-CGN is used with the preconditioner %%@
(\ref{pde2-system}) for different values of  $\alpha$.}
\begin{center}
\begin{tabular}{|c|c|c|c|c|c|c|c|}\hline
$\lambda^{s}_{max}(\frac{1-\alpha}{\alpha})$ &I& &$\lambda^{s}_{max}(\frac{1-\alpha}{\alpha})$&I %%@
\\ \hline
$\infty (\alpha=0)$& $>5000$&  &0.1 & 232 \\ \hline
$0(\alpha=1)$& 229 &  & 0.2& 237 \\ \hline
-0.1& 221&  &0.4 & 249 \\ \hline
-0.25& 216&& 0.8 &263 \\ \hline
-0.5& 201&  &1 &272\\ \hline
-0.7& 191&  &5 & 384\\ \hline
-0.8& 186&  &10 & 474 \\ \hline
-0.9& 180&  &20 & 642 \\ \hline
-0.95& 179&  &50& 890 \\ \hline
-0.99& 177&  &100& 1170 \\ \hline
-0.999& 180&  &1000& 2790 \\ \hline
-0.9999& 234&  &10000& 4807 \\ \hline

\end{tabular}
\end{center}
\end{table}
\FloatBarrier
\begin{remark}
As we see in Table 5,  for $\lambda^{s}_{max}(\frac{1-\alpha}{\alpha})=-.99$ we have the minimum %%@
number of iterations. Actually, this is the case in some other experiments, but for many other %%@
system the minimum number of iterations accrues for some other $\alpha$ with %%@
$-1<\lambda^{s}_{max}(\frac{1-\alpha}{\alpha})\leq 0$. Our experiments show that for a well chosen %%@
$\alpha>1$, one may  considerably decrease the number of iterations. Obtaining theoretical results %%@
on how to choose parameter $\alpha$ in \ref{pde2-system} seems to be an interesting problem. 
\end{remark}

Note that the coefficient matrix of the linear system corresponding to (\ref{pde2}) is positive %%@
definite. Hence we may also apply CG with the preconditioned symmetric system of equations
\begin{equation}\label{pde2-system2}
A^{T}(A_s-\alpha \lambda^{s}_{\rm min}I)^{-1}A=A^{T}(A_s-\alpha \lambda^{s}_{min}I)^{-1}b,
\end{equation} 
where $\lambda^{s}_{\rm min}$ is the smallest eigenvalue of $A_s$ and $\alpha<1$. The number of %%@
iterations function of $\alpha$, that  CG needs to converges to a solution with relative residual %%@
$10^{-6}$ are presented in Table 7.

\begin{table}[ht] \label{sol2-pde2} \caption{Number of iterations to find a solution with relative %%@
residual $10^{-6}$ for equation (\ref{pde2}) when SD-CGN is used with the preconditioner %%@
(\ref{pde2-system2}) for different values of  $\alpha$.}
\begin{center}
\begin{tabular}{|c|c|c|}\hline
$\alpha$& I  \\ \hline
0& 229\\ \hline
0.5& 204   \\ \hline
0.9& 177 \\ \hline
0.99& 166\\ \hline
0.999& 168 \\ \hline
0.9999& 181  \\ \hline
0.99999& 194  \\ \hline
0.999999& 222 \\ \hline
0.9999999& 248  \\ \hline
0.99999999& 257  \\ \hline
\end{tabular}
\end{center}
\end{table}
\FloatBarrier
\begin{remark}
As we see in the above table, for $\alpha=0.99$ in (\ref{pde2-system2}) we have the minimum number %%@
of iterations. Obtaining theoretical results on how to choose the parameter $\alpha$ seems to be %%@
an interesting problem to study.
\end{remark}
We also repeat the experiment by applying CG to the system of equations
\begin{equation}\label{pde2-system3}
A^{T}\left(A_s-0.99 \lambda^{s}_{\rm min}I)^{-1}-\frac{0.99}{\lambda^{s}_{\rm %%@
max}}I\right)A=A^{T}\left((A_s-o.99 \lambda^{s}_{\rm min}I)^{-1}-\frac{0.99}{\lambda^{s}_{\rm %%@
max}}I\right)b.
\end{equation}
Then CG needs 131 iterations to converge to a solution with relative residual $10^{-6}$.

As another experiment we apply CG to the preconditioned linear system 
\[A_s^{-1}A^{T}A_s^{-1}A=A_s^{-1}A^{T}A_s^{-1}b,\]
to solve the non-symmetric linear system obtained from discritization of the Equation %%@
(\ref{pde2}). The CG converges in 31 iterations to a solution with relative residual less than %%@
$10^{-6}$. Since, we need to solve two equations with the coefficient matrix $As$, the cost of %%@
each iteration in this case is towice as much as SD-CGN. So, by the above preconditioning we %%@
decrease cost of finding a solution to less that $62/131$ of that of SD-CGN (System %%@
(\ref{pde2-system3})). 

\begin{example} Consider now the following equation
%\begin{equation}\label{pde3}
%-\Delta u+[\partial_{x}(20\exp(3.5(x^2+y^2))u)+20\exp(3.5(x^2+y^2))\partial_{x} u]/2-200u=f(x), \ %%@
%\ \ \ \hbox{on} \ \ [0,1]\times [0,1],
%\end{equation}
\begin{equation}\label{pde3}
-\Delta u+ 10\frac{\partial (\exp(3.5(x^2+y^2)u)}{\partial x}+10\exp(3.5(x^2+y^2))\frac{\partial u}{\partial x}-200u=f(x), \ %%@
\ \ \ \hbox{on} \ \ [0,1]\times [0,1],
\end{equation}

If we discretize this equation with stepsize $1/31$ and use backward differences for the first %%@
order term,  we get a linear system of equations $Ax=b$ with $A$ being a non-symmetric and %%@
non-positive definite coefficient matrix. We then apply CG to the following preconditioned, %%@
symmetrized and positive definite matrix
\begin{equation}\label{pde3-system}
A^{T}((A_s-\alpha \lambda^{s}_{\rm min}I)^{-1}+\beta I)A=A^{T}((A_s-\alpha \lambda^{s}_{\rm %%@
min}I)^{-1}+\beta I)b,
\end{equation}
with $\alpha<1$. For different values of $\alpha$ the number of iterations which CG needs to %%@
converge to a solution with the relative residual $10^{-6}$ are presented in Table 8.
\end{example}
\begin{table}[ht] \label{sol1-pde3} \caption{Number of iterations to find a solution with relative %%@
residual $10^{-6}$ for equation (\ref{pde3}) when SD-CGN is used with the preconditioner %%@
(\ref{pde3-system}) for different values of $\alpha$ and $\beta$. }
\begin{center}
\begin{tabular}{|c|c|c|}\hline
$\alpha$& $\beta=0$&$\beta=-.99/\lambda^{s}_{max}$  \\ \hline
10& 543& 424   \\ \hline
5& 446&352 \\ \hline
2.5& 369&288 \\ \hline
1.5& 342&264 \\ \hline
1.1& 331&258 \\ \hline
1.01& 327&259\\ \hline
1.001& 333&271 \\ \hline
1.0001& 368&289  \\ \hline
1.00001& 401&317  \\ \hline
\end{tabular}
\end{center}
\end{table}
\FloatBarrier
We repeat our experiment with stepsize $1/61$ and get a system with 3600 unknowns. With %%@
$\alpha=-1.00000001$ and $\beta=0$,  CG converges  in one single iteration to a solution with %%@
relative residual less than $10^{-6}$. We also apply QMR, BiCGSTAB, BiCG, and CGS (also %%@
preconditioned  with the symmetric part as well) to solve the corresponding system of linear %%@
equations with stepsize $1/31$. The number of iterations needed to converge to a solution with %%@
relative residual $10^{-6}$ are presented in Table 9.

\begin{table}[ht] \label{sol2-pde3} \caption{Number of iterations to find a solution with relative %%@
residual $10^{-6}$ for equation (\ref{pde3}) using various algorithms.}
\begin{center}
\begin{tabular}{|c|c|c|}\hline
N=900& I \\ \hline
CGNE& $>5000$\\ \hline
QMR& 3544    \\ \hline
PQMR& 490 \\ \hline
BiCGSTAB& $>5000$\\ \hline
PBiCGSTAB& Breaks down \\ \hline
BiCG& 4527  \\ \hline
PBiCG& $>1000$  \\ \hline
CGS& 1915 \\ \hline
PCGS& 649 \\ \hline
\end{tabular}
\end{center}
\end{table}
\FloatBarrier
{\bf Acknowledgments:} This paper wouldn't have seen the light without the gentle prodding and constant encouragement of Anthony Peirce,  and  the expert guidance and generous support of Chen Greif.  They have our deep and sincere gratitude.

 \end{document}